\documentclass[12pt]{article}
\usepackage{amsmath,amssymb}
\usepackage{enumerate}
\usepackage{a4wide}

\newcommand{\N}{{\mathbb N}}
\newcommand{\R}{{\mathbb R}}

\newcommand{\Z}{{\mathbb Z}}

\newcommand{\B}{\mathcal B}

\newcommand{\A}{\mathcal A}
\newcommand{\M}{\mathcal M}

\newcommand{\cL}{{\mathcal L}}
\renewcommand{\L}{L}

\newtheorem{lemma}{Lemma}[section]
\newtheorem{proposition}[lemma]{Proposition}

\newtheorem{theorem}{Theorem}
\newcommand{\proof}{{\em Proof.}}
\newcommand{\cqfd}{\hfill$\Box$}
\newcommand{\dif}{{\rm d}}
\newcommand{\esp}{{\mathbb E}}
\newcommand{\ind}{{\mathchoice{\mathrm {1\mskip-4.1mu l}} {\mathrm{ 1\mskip-4.1mu l}} {\mathrm {1\mskip-4.6mu l}} {\mathrm {1 \mskip-5.2mu l}}}}

\newcommand{\I}{\mathcal I}

\begin{document}
\title{On Sums of Indicator Functions in\\ Dynamical Systems}
\author{Olivier Durieu\thanks{e-mail: olivier.durieu@univ-rouen.fr}
 \, \& \, Dalibor Voln\'y\thanks{e-mail: dalibor.volny@univ-rouen.fr}\\
\\
Laboratoire de Math\'ematiques Rapha\"el Salem,\\ UMR 6085 CNRS-Universit\'e de Rouen} 
\maketitle

\begin{abstract}
In this paper, we are interested in the limit theorem question for sums of indicator functions.
We show that in every invertible ergodic dynamical system, for every increasing sequence $(a_n)_{n\in\N}\subset\R_+$
such that $a_n\nearrow\infty$ and $\frac{a_n}{n}\rightarrow 0$ as $n\rightarrow\infty$, there exists a dense $G_\delta$ of measurable sets $A$ such that the sequence of the distributions of the partial sums $\frac{1}{a_n}\sum_{i=0}^{n-1}(\ind_A-\mu(A))\circ T^i$ is dense in the set of the 
probability measures on $\R$.

\medskip


\medskip

\end{abstract}

\section{Introduction}
In \cite{BurDen87}, Burton and Denker showed that in every aperiodic dynamical system
there exists a process $(f\circ T^i)$ for which the CLT holds and they posed
the question how big is the subset of $f\in \L^2$ with this property. Clearly, we
have to study the space $\L_0^2$ of $f$ with $E(f|{\mathcal I}) =0$. As already observed
by Burton and Denker, because the coboundaries are dense in $\L_0^2$, this set is
dense. In Voln\'y \cite{Vol90} it has been proved that for any sequence $a_n\to \infty$, $\frac{a_n}{n}\to 0$,
there exists a dense $G_\delta$ part $G$ of $\L_0^2$ such that for any $f\in G$
and any probability law $\nu$ there exists a sequence $n_k\to \infty$ such that
$\frac{1}{a_{n_k}}S_{n_k}(f)$ converge in law to $\nu$. The same result takes place for
all spaces $\L^p$, $1\leq p\leq \infty$. 
Liardet and Voln\'y \cite{LiaVol97} obtained the same result for the space of continuous functions for a uniquely ergodic continuous homeomorphism of a metrizable compact space and similar results for spaces of smooth
functions for irrational rotations of a circle.
As a corollary we get that generically, the rate of convergence in the ergodic
theorems (of Birkhoff and of von Neumann) may be arbitrarily slow. This gave a new
proof of a result of A. del Junco and J. Rosenblatt \cite{JunRos79}. In the paper of del Junco and
Rosenblatt a similar result on the rate of convergence in the ergodic theorems
was found for functions $\ind_A - \mu(A)$, the genericity was studied in the
space of $A\in\A$ equipped with the (pseudo)metric of the measure of
symmetric difference.

In the present paper we shall study the distributional convergence for the
functions $\ind_A - \mu(A)$. The research was motivated by the study of the
invariance principle of the empirical process of strictly stationary sequences
$(X_i)_{i\in\N}$ in Dehling, Durieu and Voln\'y \cite{DehDurVol08}.

\medskip

\section{Result}

Let $(\Omega,\A,\mu)$ be a non-atomic Lebesgue probability space and $T$ be an invertible measurable transformation from $\Omega$ to $\Omega$. We say that $T$ is measure preserving if
for all $A\in\A$, $\mu(T^{-1}A)=\mu(A)$. 
In the sequel, we will often say that $(\Omega,\A,\mu,T)$ is a dynamical system when $(\Omega,\A,\mu)$ is a non-atomic Lebesgue probability space and $T$ is an invertible measure preserving transformation of $\Omega$.

\medskip

Further, the transformation $T$ is ergodic if $T^{-1}A=A$ implies that $\mu(A)=0$ or $1$. It is aperiodic if 
$$
\mu\{x\in\Omega\,/\,\exists n\ge 1, T^{n}x=x\}=0.
$$

\medskip

On $\A$ we consider the pseudo-metric $\Theta$ defined by
$$
\Theta(A,B)=\mu(A\triangle B),\quad   A,B\in \A.
$$
Our main result is the following, where $X_n\xrightarrow[\;n\to\infty\;]{\mathcal D}\nu$ means that the sequence of real random variables $X_n$ converges in distribution to a real random variable having distribution $\nu$.

\begin{theorem}\label{M2}
 Let $(\Omega,\A,\mu)$ be a non-atomic Lebesgue probability space and $T$ be an ergodic invertible measure preserving transformation of $\Omega$.
Let $(a_n)_{n\in\N}\subset\R_+$ be an increasing sequence
such that $a_n\nearrow\infty$ and $\frac{a_n}{n}\rightarrow 0$ as $n\rightarrow\infty$.

There exists a $\Theta$-dense $G_\delta$ of sets $A\in\A$ having the property that for every probability $\nu$ on $\R$, there exists a subsequence $(n_k)_{k\in\N}$
satisfying
\begin{equation*}
 \frac{1}{a_{n_k}}S_{n_k}(\ind_A-\mu(A))\;\xrightarrow[\;k\to\infty\;]{\mathcal D}\;\nu.
\end{equation*}
\end{theorem}

To be complete, notice that, in the non-ergodic case, to find a set which satisfies the conclusion of Theorem \ref{M2}, we have to consider the sets $A$ such that $\esp(\ind_A|\I)$ is almost surely constant. The class of such sets is not, in general, dense in $\A$. So, in the non-ergodic case, we cannot expect the result of genericity.

Nevertheless, we can prove the existence of such sets by an explict construction. This is the purpose of the paper by Durieu and Voln\'y \cite{DurVol08c}.

\medskip

\section{Some preliminary results}

Let $\M$ be the set of all probability measures on $\R$ and
$\M_0$ be the set of all probability measures on $\R$ which have zero-mean.
Recall that $\M_0$ is dense in $\M$ for the topology of the weak convergence. 
We denote by $d$ the L\'evy metric on $\M$. For all $\mu$ and $\nu$ in $\M$ with distribution functions
$F$ and $G$,
$$
d(\mu,\nu)=\inf\{\varepsilon>0\,/\,G(t-\varepsilon)-\varepsilon\le F(t)\le G(t+\varepsilon)+\varepsilon, \forall t\in\R\}.
$$
The space $(\M,d)$ is a complete separable metric space and convergence with respect to $d$ is equivalent to weak convergence of distributions (see Dudley \cite{Dud02}, pages 394-395).

\medskip

If $X:\Omega\longrightarrow\R$ is a random variable, we denote by $\cL_\Omega(X)$ the distribution of $X$ on $\R$.

\begin{lemma}\label{function}
 Let $(\Omega,\A,\mu)$ be a Lebesgue probability space and $\nu$ be a probability on $\R$.
Then, there exists a random variable $X:\Omega\longrightarrow\R$, such that
$$
\cL_{\Omega}(X)=\nu.
$$
\end{lemma}

\proof

It is well known that  $(\Omega,\A,\mu)$ is isomorphic to $([0,1],\B[0,1],\lambda)$, where $\lambda$ is the Lebesgue measure on $[0,1]$. If $Q_\nu$ denotes the pseudo-inverse of the distribution function of $\nu$ and $U$ is the identity on $[0,1]$, it is classical that $\cL_{[0,1]}(Q_\nu(U))=\nu$.
\cqfd

\medskip

Let $\nu$ be a probability on $\R$. For $B\in\B(\R)$ with $\nu(B)>0$, $\nu_B$ denotes the probability on $\R$ defined by
$$
\nu_B(A)=\nu(B)^{-1}\nu(A\cap B).
$$
For $x\in\R$, $\nu_x$ denotes the probability on $\R$ defined by
$$
\nu_x(B)=\nu(xB)
$$
where $xB=\{xb\,/\,b\in B\}$.

Here are some properties of the L\'evy metric which will be used in the sequel.

\begin{lemma}\label{Lnot}
$\,$
\begin{enumerate}[{\rm (i)}]
 \item For each probability $\nu$ on $\R$, for all Borel sets $B$,
$$
d(\nu_B,\nu)\le \nu(\R\smallsetminus B).
$$
\item For all probabilities $\nu$ and $\eta$ on $\R$, for all $x\ge 1$,
$$
d(\nu_x,\eta_x)\le d(\nu,\eta).
$$
\item For all probability $\nu$ on $\R$, for all mesurable functions $f$ and $g$ from $\Omega$ to $\R$,
$$
d(\cL_\Omega(f+g),\nu)\le (\cL_\Omega(f),\nu)+d(\cL_\Omega(g),\delta_0)
$$
where $\delta_0$ is the Dirac measure at $0$.
\item For all probability $\nu$ on $\R$, 
$$
d(\nu,\delta_0)\le A \mbox{ if and only if } \nu((-\infty,-A))\le A \mbox{ and } \nu((A,\infty))\le A.
$$
\end{enumerate}
\end{lemma}

The proof is an exercise which is left to the reader.

\medskip

\begin{lemma}\label{Lsupport}
Let $(a_n)_{n\in\N}\subset\R$ be an increasing sequence
such that $a_n\nearrow\infty$ as $n\rightarrow\infty$.
 For each probability $\nu\in\M_0$ and $\varepsilon>0$, there exist $C\ge 1$ and $n_0\in\N$, for all $n\ge n_0$, there exists a probability $\eta$ on $\R$ with support $S\subset [-a_nC,a_nC]\cap\Z$ such that
$$
d(\eta_{a_n},\nu)\le \varepsilon
$$
and
$$
\esp(\eta):=\int x{\rm d}\eta=0.
$$
\end{lemma}

\proof

Let $\nu\in\M_0$ and $\varepsilon>0$ be fixed and choose $\alpha>0$ such that $6\alpha\le\varepsilon$ and $\alpha<\frac{1}{2}$.
There exists $C\ge 1$ such that 
$$
\int_{\R\smallsetminus[-C,C]}|x|\dif\nu(x)\le \alpha. 
$$
In particular, $\nu(\R\smallsetminus[-C,C])\le\alpha$.
Define $\tau=\nu_{[-C,C]}$. By Lemma \ref{Lnot},
$$
d(\tau,\nu)\le \nu(\R\smallsetminus[-C,C])\le\alpha
$$
and we have
$$
|\esp(\tau)|\le \nu([-C,C])^{-1}|\esp(\nu)-\int_{\R\smallsetminus[-C,C]}x\dif\nu(x)|\leq \frac{\alpha}{1-\alpha}.
$$

\medskip

Now, choose $n_0\in\N$ such that $\frac{1}{a_{n_0}}<\alpha$ and fix $n\ge n_0$.
Then we define the probability $\eta'$ on $\R$ with support in $\Z$, by $\eta'(\{k\}):=\tau\left(\left[\frac{k}{a_n},\frac{k+1}{a_n}\right)\right)$, $k\in\Z$.

We have, for all $t\in\R$,
\begin{eqnarray*}
 \eta'_{a_n}((-\infty,t])&=&\eta'((-\infty,\lfloor ta_n\rfloor])\\
&=&\tau\left(\left(-\infty,\frac{\lfloor ta_n\rfloor+1}{a_n}\right)\right)\\
&\le& \tau\left(\left(-\infty,t+\frac{1}{a_n}\right]\right)
\end{eqnarray*}
and
$$
\tau((-\infty,t])\le\eta'_{a_n}((-\infty,t]).
$$
Thus $d(\eta'_{a_n},\tau)\le \frac{1}{a_n}\le \alpha$ and
$d(\eta'_{a_n},\nu)\le 3\alpha$.

\medskip

So, if $\esp(\eta')=0$, $\eta'$ verifies the conclusion of the proposition. If it is not the case, we proceed as follows.
Observe that
$$
a_n\esp(\tau)-1\le \esp(\eta')\le a_n\esp(\tau),
$$
and thus, since $\alpha<\frac{1}{2}$ and $a_n\alpha>1$,
$$
|\esp(\eta')|\le a_n|\esp(\tau)|+1\le a_n\frac{\alpha}{1-\alpha}+1\le 3 a_n\alpha.
$$
We denote by $s$ the sign of $\esp(\eta')$ and we set $p=1+\frac{|\esp(\eta')|}{\lfloor a_nC\rfloor}$.

\medskip

Now we denote by $\eta$ the probability on $\R$ with support in $\{-\lfloor a_n C\rfloor,\dots,\lfloor a_nC\rfloor\}$ defined by
$$
\eta(\{i\})=\left\{
\begin{array}{lll}
 \frac{1}{p}(\eta'(\{i\})+\frac{|\esp(\eta')|}{\lfloor a_nC\rfloor})&\mbox{ if } & i= -s\lfloor a_nC\rfloor \\
\frac{1}{p}\eta'({i})&\mbox{ otherwise }
\end{array}
\right.
$$
Then $\esp(\eta)=0$ and by Lemma \ref{Lnot} (ii),
$$
d(\eta'_{a_n},\eta_{a_n})\le d(\eta',\eta)\leq \frac{|\esp(\eta')|}{\lfloor a_nC\rfloor}\le 3\alpha.
$$
Therefore $d(\eta_{a_n},\nu)\le \varepsilon$.
\cqfd

\medskip

\begin{proposition}\label{lem}
Let $(\Omega,\A,\mu,T)$ be an ergodic dynamical system, $(a_n)_{n\in\N}\subset\R_+$ be an increasing sequence
such that $a_n\nearrow\infty$ and $\frac{a_n}{n}\rightarrow 0$ as $n\rightarrow\infty$ and $\varepsilon>0$.
Let $A\in\A$ be a set such that $\mu(A)<1$ and $\nu$ be a probability in $\M_0$.
There exists $N\in\N$ such that for any $n\ge N$, there exists a set $B_n\in\A$ such that $\mu(B_n)\le \varepsilon$, $A\cap B_n=\emptyset$ and
$$
d(\cL_\Omega(\frac{1}{a_n}S_n(\ind_{B_n}-\mu(B_n))),\nu)\le\varepsilon
$$
\end{proposition}

\proof

Fix $\varepsilon>0$ and $A\in\A$ such that $\mu(A)<1$. Let $\alpha$ be a positive constant such that $\alpha\le\frac{\varepsilon}{5}$ and $\mu(A)+2\alpha<1$.

\medskip

By Lemma \ref{Lsupport} applied to $\nu$ and $\alpha$ we get the constants $C\ge 1$ and $n_0\ge 1$ for which the conclusion of the lemma holds. Set $\gamma:=\frac{\alpha}{C+1}$.
Let $n_1$ be an integer such that, for all $n\ge n_1$,
\begin{equation}\label{equ1}
2\frac{a_nC+1}{n}\le \alpha.
\end{equation}

\medskip

Applying Birkhoff's ergodic theorem and Egorov's theorem, we get that there exist a set $E\in\A$ of measure greater than $1-\frac{\gamma}{2}$ and an integer $n_2\ge 1$ such that for all
$n\ge n_2$, for all $x\in E$,
\begin{equation}\label{equ2}
\left|\frac{1}{n}S_n(\ind_A-\mu(A))(x)\right|\le \alpha.
\end{equation}

\medskip

We denote by $\overline{n}$ the maximum of $n_0$, $n_1$ and $n_2$ and we choose $N\in\N$ such that
$$
\frac{\overline{n}}{N}\le \alpha.
$$

\medskip

For any $n\ge N$, there exists a Rokhlin tower of height $n$ with base $F\subset E$ and junk set of measure smaller than $\gamma$.

\medskip

Indeed, let $G$ be the base of a Rokhlin tower of height $n$ and of measure greater than $1-\frac{\gamma}{2}$.
Because $\mu(\Omega\smallsetminus E)\le\frac{\gamma}{2}$, there exists an integer $i_0\in\{0,\dots,n-1\}$ such that 
$$
\mu((T^{i_0}G)\cap E)\ge \frac{1}{n}(1-\frac{\gamma}{2}-\frac{\gamma}{2}))= \frac{1-\gamma}{n}.
$$
If $F=T^{i_0}G\cap E$, then $F\subset E$ and the sets $F$, $TF$,...,$T^{n-1}F$ are disjoint. So, $F$ is the base of a Rokhlin tower of height $n$ with a junk set of measure
smaller than $\gamma$.

\medskip

From now on, $n$ is fixed.
By Lemma \ref{Lsupport}, there exists a centered probability $\eta$ with support in $S\subset [-a_nC,a_nC]\cap \Z$ such that
$$
d(\eta_{a_n},\nu)\le \alpha.
$$
By Lemma \ref{function}, there exists a function $h:F\longrightarrow \Z$ such that
$\cL_F(h)=\eta$. In particular, $\esp_F(h):=\int_F h {\rm d}\mu_F=0$.
We set $d=\lfloor a_nC\rfloor +1$ and
$$
\begin{array}{rcl}
 g:F &\longrightarrow &\Z\\
 x &\longmapsto & h(x)+d.
\end{array}
$$
Note that $1\le g\le 2d$ almost surely and $\esp_F(g)=d$.
We now set 
$$
F_i:=g^{-1}(\{i\}),\; i=1,\dots,2d.
$$
Note that the $F_i$'s depend on $\nu$, $\alpha$, $C$ and $n$.
Further
$\{F_1, F_2,\dots,F_{2d}\}$ is a partition of the set $F$.

\medskip

By (\ref{equ2}), for each $x\in F$, the sub-orbit $\{x,Tx,\dots,T^{\overline{n}-1}x\}$ hits $A$ at most
$\overline{n}(\mu(A)+\alpha)$ times (\textit{i.e.} $\#(\{x,Tx,\dots,T^{\overline{n}-1}x\}\cap A)\le \overline{n}(\mu(A)+\alpha)$. Since by (\ref{equ1}), $\frac{2d}{\overline{n}}\le \alpha$ and $\mu(A)+2\alpha<1$, 
we can find $2d$ points in this sub-orbit which are not in $A$.

Then, for each $i=1,\dots,2d$ and for each $x\in F$, we can denfine the set $b_i(x)$ composed by the $i$ first points of the sub-orbit $\{x,Tx,\dots,T^{\overline{n}-1}x\}$ which are not in $A$.

We now set, for each $i=1,\dots,2d$, 
$$
B_i=\bigcup_{x\in F_i}b_i(x).
$$
Thus the $B_i$'s are disjoint measurable sets.  For each $i=1,\dots,2d$, $B_i\subset\bigcup_{j=0}^{\overline{n}-1}T^jF_i$,
$B_i\cap A=\emptyset$, $\mu(B_i)=i\mu(F_i)$ and for any $x\in F_i$, 
$$
S_{\overline{n}}(\ind_{B_i})(x)=i.
$$

\medskip

Finally, we set
$$
B=\bigcup_{i=1}^{2d}B_i.
$$
We have $B\in\A$ and $A\cap B=\emptyset$.

\medskip

From the construction of $B$ and (\ref{equ1}),
$$
 \mu(B)=\esp_{F}(g)\mu(F)= d\mu(F)\le\frac{a_{n}C+1}{n}\le\alpha\le\varepsilon.
$$

\medskip

We define 
$$
\Omega_{k}=\bigcup_{i=0}^{n-\overline{n}-1}T^{-i}F_k,\quad k=1,\dots,2d
$$
and
$$
\overline{\Omega}=\bigcup_{k=1}^{2d}\Omega_{k}.
$$
Since the $T^{-i}F_k$ are disjoint, using (\ref{equ1}) and the fact that $\gamma\le\alpha$, we have
\begin{equation}\label{omega}
\mu(\overline{\Omega})=(n-\overline{n})\mu(F)\ge 1-\gamma-\frac{\overline{n}}{n}\ge 1-3\alpha.
\end{equation}

\medskip

For $x\in\overline{\Omega}$, by construction and by disjointness of the $T^{-i}F_k$, we have
$$
S_{n}(\ind_{B})(x)=k \mbox{ if and only if } x\in\Omega_{k}.
$$
Therefore, for all $k=1,\dots,2d$,
$$
\mu_{\overline{\Omega}}(S_{n}(\ind_{B})=k)=\mu_{\overline{\Omega}}(\Omega_{k})=\mu_{F}(F_k).
$$
Thus $\cL_{\overline{\Omega}}(S_{n}(\ind_{B}))=\cL_F(g)$ and by (\ref{omega}) and Lemma \ref{Lnot} (i),
$$
d(\cL_{\Omega}(S_{n}(\ind_{B})),\cL_F(g))\le 3\alpha.
$$
So, by Lemma \ref{Lnot} (ii),
$$
d(\cL_{\Omega}(\frac{1}{a_{n}}S_{n}(\ind_{B})),\cL_F(\frac{g}{a_{n}}))\le 3\alpha
$$
and
\begin{equation}\label{e1}
d(\cL_{\Omega}(\frac{1}{a_{n}}(S_{n}(\ind_{B})-d)),\cL_F(\frac{g-d}{a_{n}}))\le 3\alpha.
\end{equation}
Now, remark that
\begin{equation}\label{e2}
d(\cL_{\Omega}(\frac{1}{a_{n}}(S_{n}(\ind_{B})-d)),\cL_{\Omega}(\frac{1}{a_{n}}S_{n}(\ind_{B}-\mu(B))))\le \alpha.
\end{equation}
Indeed, since $\mu(B)=d\mu(F)$, we have
$$
-\gamma d\le n\mu(B)-d\le 0
$$
and then
$$
\left|\frac{1}{a_{n}}S_{n}(\ind_{B}-\mu(B))-
\frac{1}{a_{n}}(S_{n}(\ind_{B})-d)\right|\le \frac{\gamma d}{a_{n}}
\le\alpha.
$$

\medskip

To conclude, using (\ref{e1}), (\ref{e2}) and the fact that $d(\cL_F(\frac{h}{a_{n}}),\nu)\le \alpha$, we get
$$
d(\cL_{\Omega}(\frac{1}{a_{n}}S_{n}(\ind_{B}-\mu(B))),\nu)\le 5\alpha\le \varepsilon.
$$
\cqfd

\medskip

\begin{proposition}\label{thelemma}
 Let $(\Omega,\A,\mu,T)$ be an ergodic dynamical system, $(a_n)_{n\in\N}\subset\R_+$ be an increasing sequence
such that $a_n\nearrow\infty$ and $\frac{a_n}{n}\rightarrow 0$ as $n\rightarrow\infty$ and $\varepsilon>0$.
For any set $A\in\A$ such that $\mu(A)<1$,
there exists a set $B\in\A$ such that 
\begin{enumerate}[{\rm (i)}]
\item $\mu(A\triangle B)\le\varepsilon$,
\item there exists a sequence $(n_k)_{k\ge 1}$ such that for all $k\ge 1$,
$$
d(\cL_{\Omega}(\frac{1}{a_{n_k}}S_{n_k}(\ind_{B}-\mu(B))),\delta_0)\le \varepsilon.
$$
\end{enumerate}
\end{proposition}

This proposition will be proved as a corollary of the following lemma.

\begin{lemma}\label{correct}
Let $(\Omega,\A,\mu,T)$ be an ergodic dynamical system, $(a_n)_{n\in\N}\subset\R_+$ be an increasing sequence
such that $a_n\nearrow\infty$ and $\frac{a_n}{n}\rightarrow 0$ as $n\rightarrow\infty$ and $\varepsilon>0$.
For any set $A\in\A$ such that $\mu(A)<1$, for any $N\in\N$, there exist $n\ge N$ and a set $C\in\A$ such that
$\mu(A\triangle C)\le\varepsilon$ and
$$
d(\cL_{\Omega}(\frac{1}{a_{n}}S_{n}(\ind_{C}-\mu(C))),\delta_0)\le \varepsilon.
$$
\end{lemma}

\proof

Let $\varepsilon>0$, $A\in\A$ such that $\mu(A)<1$ and $N\in\N$ be fixed and let $\alpha$ be a positive constant such that
 $5\alpha\le \varepsilon$ and $\mu(A)+2\alpha<1$.

\medskip

By Birkhoff's ergodic theorem and Egorov's theorem, there exist a set $G\in\A$ with $\mu(G)>1-\alpha$ and an integer $M$ such that 
for all $k\ge M$ and for all $x\in G$,
\begin{equation}\label{A}
\left|\frac{1}{k}S_k(\ind_A-\mu(A))(x)\right|\le\alpha.
\end{equation}
Furthermore, we can choose $M$ such that $\frac{1-\alpha}{M}\le \alpha$.

\medskip
There exists an integer $n\ge N$, such that
\begin{equation}\label{ttttt}
\frac{M}{a_n}\le \alpha.
\end{equation}

\medskip

 Let $F$ be the base of a Rokhlin tower of height $Mn$ with a junk set of measure smaller than $\alpha$.
For $k=0,\dots,M-1$, we define the sets 
$$
F_k=\bigcup_{i=0}^{n-1}T^{iM+k}F.
$$
Since $\mu(G)> 1-\alpha$, there exists $k_0\in\{0,\dots,M-1\}$ such that
$$
\mu(F_{k_0}\smallsetminus G)\le \frac{2\alpha}{M}.
$$
Further, $H=T^{k_0}F$ is the base of a Rokhlin tower of height $Mn$ with a junk set $J$ such that $\mu(J)\le \alpha$.

\medskip

 For $x\in \Omega$, we denote by $s_l(x)=\{x,Tx,\dots,T^{l-1}x\}$ the sub-orbit of length $l$ which begins at $x$. The set $A$ will be modified in a way that for each orbit $s_M(x)$, $x\in F_{k_0}$, the average of visits of the set along the orbit be close to the measure $\mu(A)$.

\medskip

For $x\in H$, we modify the set $A$ along the
sub-orbit $s_M(x)$ in the following way. We write $a(x):=A\cap s_M(x)$. There are two situations:

\medskip

$\bullet$ If $x\in G$, by (\ref{A}), along the sub-orbit $s_M(x)$, the number of visits of the set $A$ belongs to
$[M(\mu(A)-\alpha),M(\mu(A)+\alpha)]$ (\textit{i.e.} $\# a(x)\in [M(\mu(A)-\alpha),M(\mu(A)+\alpha)]$).
Further, $M(\mu(A)+\alpha)\le (1-\alpha)M$. So by adding or removing at most $\alpha M$ points to $a(x)$, we can modify the set $a(x)$ to get a set $a_0(x)$ (with $a_0(x)\subset s_M(x)$) such that along the sub-orbit $s_M(x)$, the number of visits of the set $a_0(x)$ belongs to $[M\mu(A)-1,M\mu(A)+1]$. A way to do that is to remove or add the $|\# a(x)-\lfloor M\mu(A)\rfloor|$ first points of $a(x)$ depending on whether $\# a(x)-\lfloor M\mu(A)\rfloor$ is positive or not.

\medskip

$\bullet$ If $x\notin G$, we can also modify $a(x)$ in order to have that the number of visits of the set $a_0(x)$ along the segment $s_M(x)$ belongs to $[M\mu(A)-1,M\mu(A)+1]$. Here, to do that we possibly need to add or remove $M$ points to $a(x)$ along $s_M(x)$.

\bigskip

To summarize, for each $x\in H$, we can modify the set $A\cap s_M(x)$ to get a set $a_0(x)\subset s_M(x)$ such that 
$\# a_0(x)\in [M\mu(A)-1,M\mu(A)+1]$ and if $x\in G$, $\# (a_0(x)\triangle a(x))\le \alpha M$, if $x\notin G$, $\# (a_0(x)\triangle a(x))\le  M$.

\medskip

We then have a set $A_0=\bigcup_{x\in H} a_0(x)$  having the property that for all $x\in H$,
\begin{equation}\label{qsz}
\left|S_M(\ind_{A_0}-\mu(A))(x)\right|\le 1
\end{equation}
(the problem of measurability for $A_0$ is not discussed, but we can see that the modifications can be done in a measurable way).

\bigskip

Now we will do almost the same modifications on sub-orbits of length $M$ starting from $T^MH$. 

For each $x\in T^MH$, we modify the set $a(x)=A\cap s_M(x)$ to get a set $a_1(x)\subset s_M(x)$ such that $\# a_1(x)\in [M\mu(A)-1,M\mu(A)+1]$. 

Notice that for each $x\in T^MH$, considering the number $\# a_0(T^{-M}x)$, we can further define $a_1(x)$ in such a way that $\# (a_0(T^{-M}x)\cup a_1(x))\in [2M\mu(A)-1,2M\mu(A)+1]$.
We keep the property that if $x\in T^MH\cap G$, $\# (a_1(x)\triangle a(x))\le \alpha M$, if $x\in T^MH\smallsetminus G$, $\# (a_1(x)\triangle a(x))\le  M$. 

We can define $A_1=\bigcup_{x\in T^MH}a_1(x)$. Then for all $x\in T^MH$,
\begin{equation}
\left|S_M(\ind_{A_1}-\mu(A))(x)\right|\le 1
\end{equation}
and for all $x\in H$,
\begin{equation*}
 \left|S_{2M}(\ind_{A_0\cup A_1}-\mu(A))(x)\right|\le 1.
\end{equation*}

\bigskip

Now we do the same modifications for all points of $T^{2M}H$, $T^{3M}H$, ... , and $T^{(n-1)M}H$. Finally, we can get a measurable set $B=\bigcup_{i=0}^{n-1}A_i$ deduced from $A$ and having the property that for all $k\in \{1,\dots,n\}$ and for each $x\in T^{kM}H$,
\begin{equation}
\left|S_M(\ind_{B}-\mu(A))(x)\right|\le 1
\end{equation}
Further, for all $x\in H$ and for all $k\in \{1,\dots,n\}$,
\begin{equation*}
 \left|S_{kM}(\ind_{B}-\mu(A))(x)\right|\le 1.
\end{equation*}
Note that we did not change the set $A$ on the junk set $J$, so $A\cap J=B\cap J$. Recall that we did at most $\alpha M$ changes for points
of $F_{k_0}\cap G$ and at most $M$ for points of $F_{k_0}\smallsetminus G$, so
\begin{eqnarray}
\mu(A\triangle B)
&\le & \alpha M \mu(F_{k_0}\cap G)+ M\mu(F_{k_0}\smallsetminus G)\nonumber\\
&\le & 3\alpha. \label{jk}
\end{eqnarray}
We also have
\begin{eqnarray*}
|\mu(B)-\mu(A)|
&\le & |\int_{H}S_{Mn}(\ind_{B})(x)\dif \mu(x)+\mu(B\cap J)-\mu(A)Mn\mu(H)-\mu(A)\mu(J)|\\
&\le& \int_{H}|S_{Mn}(\ind_{B}-\mu(A))(x)|\dif \mu(x)+|\mu(B\cap J)-\mu(A)\mu(J)|\\
&\le& \mu(H) + |\mu(B\cap J)-\mu(A)\mu(J)|\\
&\le& \mu(H) + \mu(J).
\end{eqnarray*}
Then, by changing the set $B$ on only one level of the tower and on the junk set,
we can obtain a new set $C\in\A$ such that $\mu(A)=\mu(C)$ and $\mu(B\triangle C)\le \mu(H) + \mu(J) \le 2\alpha$.

Thus, by (\ref{jk}), 
$$
\mu(A\triangle C)\le 5 \alpha \le \varepsilon
$$
and we have the following property:
for all $x\in H$, for all $k,k'\in\{1,\dots,n\}$,
\begin{equation}\label{yyyyy}
 \left|S_{kM}(\ind_{C}-\mu(C))(x)-S_{k'M}(\ind_{C}-\mu(C))(x)\right|\le 3.
\end{equation}

\medskip

Let $\tilde{\Omega}=\bigcup_{i=0}^{n(M-1)}T^iH$. We have 
$$\mu(\Omega\smallsetminus\tilde{\Omega})\le n\mu(H) + \mu(J)\le 2\alpha$$
and, by (\ref{yyyyy}) and (\ref{ttttt}), for all $x\in\tilde{\Omega}$,
$$
\frac{1}{a_n}\left|S_{n}(\ind_{C}-\mu(C))(x)\right|\le \frac{3 +2M}{a_n}\le 3\alpha
$$
From these two inequalities, by Lemma \ref{Lnot} (iv), we deduce that
$$
d(\cL_{\Omega}(\frac{1}{a_{n}}S_{n}(\ind_{C}-\mu(C))),\delta_0)\le 3\alpha\le \varepsilon.
$$
\cqfd

\medskip

{\em Proof of Proposition \ref{thelemma}.}

Let $\varepsilon>0$ and $A\in\A$ such that $\mu(A)<1$ be fixed and choose $\varepsilon_1\le \frac{\varepsilon}{2}$.

\medskip

By Lemma \ref{correct}, there exist $n_1\in\N$ and a set $C_1\in\A$ such that $\mu(A\triangle C_1)\le \varepsilon_1$ and
$$
d(\cL_{\Omega}(\frac{1}{a_{n_1}}S_{n_1}(\ind_{C_1}-\mu(C_1))),\delta_0)\le \varepsilon_1.
$$

\medskip

We will proceed by induction. After step $k-1$, we choose $\varepsilon_{k}= \frac{\varepsilon_{k-1}}{2n_{k-1}}$.
By application of Lemma \ref{correct}, there exist an integer $n_{k}\ge n_{k-1}$ and a set $C_{k}\in\A$ such that $\mu(C_{k-1}\triangle C_{k})\le \varepsilon_{k}$ and
$$
d(\cL_{\Omega}(\frac{1}{a_{n_k}}S_{n_k}(\ind_{C_k}-\mu(C_k))),\delta_0)\le \varepsilon_k.
$$
Note that for all $i,j>0$, $\varepsilon_{i+j}\le\frac{\varepsilon_i}{2^j}$.

\medskip

Finally, we set $B=\bigcup_{n\ge 1}\bigcap_{k\ge n} C_k$. Noticing that $D\triangle (E\cap F)\subset (D\triangle E)\cup (E\triangle F)$ for all sets D,E,F, we get
$$
\mu(A\triangle B)\le \sum_{k\ge 1}\varepsilon_k\le\varepsilon
$$
and for all $k\ge 1$, using Lemma \ref{Lnot} (iii) and (iv),
\begin{eqnarray*}
d(\cL_{\Omega}(\frac{1}{a_{n_k}}S_{n_k}(\ind_{B}-\mu(B))),\delta_0)
&\le & d(\cL_{\Omega}(\frac{1}{a_{n_k}}S_{n_k}(\ind_{C_k}-\mu(C_k))),\delta_0)\\
&&+ d(\cL_{\Omega}(\frac{1}{a_{n_k}}S_{n_k}(\ind_{B\smallsetminus C_k}-\mu(B\smallsetminus C_k))),\delta_0)\\
&&+ d(\cL_{\Omega}(\frac{1}{a_{n_k}}S_{n_k}(\ind_{C_k\smallsetminus B}-\mu(C_k\smallsetminus B))),\delta_0)\\
&\le& \varepsilon_k+ n_k\mu(B\triangle C_k)\\
&\le& \varepsilon_k+ n_k\sum_{i\ge k+1}\mu(C_{i-1}\triangle C_{i})\\
&\le& \varepsilon_k+ n_k\sum_{i\ge k+1}\varepsilon_i\\
&\le& \varepsilon.
\end{eqnarray*}
\cqfd

\medskip

\section{Proof of Theorem \ref{M2}}\label{STHEO}

Let $(\Omega,\A,\mu,T)$ be an ergodic dynamical system. Let $(a_n)_{n\in\N}\subset\R_+$ be an increasing sequence
such that $a_n\nearrow\infty$ and $\frac{a_n}{n}\rightarrow 0$ as $n\rightarrow\infty$.

Let $(\varepsilon_k)_{k\ge 1}$ be a decreasing sequence of positive reals such that $\varepsilon_k$ goes to $0$ as $k$ goes to $\infty$.

\medskip

For each $\nu\in\M_0$ and for each $k\ge 1$, we define
$$
H_k^\nu=\{A\in\A\,/\,\exists n\ge k\mbox{ such that } 
d(\cL_\Omega(\frac{1}{a_{n}}S_{n}(\ind_{A}-\mu(A))),\nu)< \varepsilon_k\}.
$$
For each $\nu\in\M_0$ and for each $k\ge 1$, it is clear that $H_k^\nu$ is an open set in $\A$.
We now prove that it is dense. 

\medskip

Assume that $\nu$ and $k$ are fixed and let $\varepsilon>0$ and $A\in\A$.
By Proposition \ref{thelemma}, there exists a set $B\in\A$ such that $\mu(A\triangle B)<\frac{\varepsilon}{2}$ and 
there exists a sequence $(n_i)_{i\ge 1}$ such that for all $i\ge 1$,
\begin{equation}\label{zzz1}
d(\cL_{\Omega}(\frac{1}{a_{n_i}}S_{n_i}(\ind_{B}-\mu(B))),\delta_0)\le \frac{\varepsilon_k}{2}.
\end{equation}

By Proposition \ref{lem},
there exists an integer $i_0$ such that, for the integer $n=n_{i_0}\ge k$, there exists a set $C\in\A$ satisfying $\mu(C)<\frac{\varepsilon}{2}$,
$C\cap B=\emptyset$ and
\begin{equation}\label{zzz2}
d(\cL_\Omega(\frac{1}{a_{n}}S_{n}(\ind_{C}-\mu(C))),\nu)< \frac{\varepsilon_k}{2}.
\end{equation}

\medskip

Hence, $\mu((B\cup C) \triangle A)<\varepsilon$ and since $B$ and $C$ are disjoint, by Lemma \ref{Lnot} (iii) and by (\ref{zzz1}) and (\ref{zzz2}), we get
\begin{eqnarray*}
&&\hspace{-30pt} d(\cL_\Omega(\frac{1}{a_{n}}S_{n}(\ind_{B\cup C}-\mu(B\cup C))),\nu)\\
&&\le d(\cL_\Omega(\frac{1}{a_{n}}S_{n}(\ind_{C}-\mu(C))),\nu)
+d(\cL_\Omega(\frac{1}{a_{n}}S_{n}(\ind_{B}-\mu(B))),\delta_0)\\
&&< \varepsilon_k,
\end{eqnarray*}
\textit{i.e.} $B\cup C$ belongs to $H_k^\nu$.
Therefore $H_k^\nu$ is dense in $\A$ for the pseudo-metric $\Theta$. 

\medskip

Let M be a countable subset of $\M_0$ which is dense in $\M$ and set 
$$
H=\bigcap_{\nu\in M}\bigcap_{k=1}^{\infty}H_k^\nu.
$$
By Baire's theorem, $H$ is a dense $G_\delta$ (for the metric of the measure of the symmetric difference).

\medskip

Further, for each $A\in H$, the sequence $(\cL_\Omega(\frac{1}{a_{n}}S_{n}(\ind_A-\mu(A))))_{n\ge 1}$ is dense in $\M$ for the L\'evy metric $d$.

Indeed, let $A\in H$, $\eta\in \M$ and $\varepsilon>0$.
By density of $M$, there exist $\nu\in M$ such that 
$$
d(\nu,\eta)<\frac{\varepsilon}{2}.
$$
But $A\in H_k^\nu$ for all $k\ge 1$, then there exists an increasing sequence $(n_k)_{k\ge 1}$ such that
\begin{equation*}
 d(\cL_\Omega(\frac{1}{a_{n_k}}S_{n_k}(\ind_A-\mu(A))),\nu)\le \varepsilon_k.
\end{equation*}
Thus, there exists $K\in\N$ such that
\begin{equation*}
 d(\cL_\Omega(\frac{1}{a_{n_{K}}}S_{n_{K}}(\ind_A-\mu(A))),\nu)\le \frac{\varepsilon}{2}
\end{equation*}
and
\begin{equation*}
 d(\cL_\Omega(\frac{1}{a_{n_{K}}}S_{n_{K}}(\ind_A-\mu(A))),\eta)\le \varepsilon.
\end{equation*}
\cqfd

\section*{Acknowledgement}

We would like to thank Jean-Pierre Conze who found an error in an earlier version of this paper and S\'ebastien Gou\"ezel who suggested us a simplification in the proof which is done in Section \ref{STHEO}. We also thank the referee for his advices.

\small
\bibliographystyle{plain}

\end{document}